\documentclass{amsart}
\usepackage{amssymb,latexsym,epsfig,color}

\newtheorem{theorem}{Theorem}[section]
\newtheorem{proposition}[theorem]{Proposition}

\newtheorem{conjecture}[theorem]{Conjecture}
\newtheorem{lemma}[theorem]{Lemma}

\theoremstyle{definition}
\newtheorem{definition}[theorem]{Definition}
\newtheorem{remark}[theorem]{Remark}
\newtheorem{example}[theorem]{Example}

\def\proof{\smallskip\noindent {\it Proof: \ }}
\def\endproof{\hfill$\square$\medskip}

\def\F{\mathbb{F}}
\def\N{\mathbb{N}}
\def\M{\mathcal{M}}

\newcommand{\reg}{\mbox{\upshape reg}}
\newcommand{\lk}{\mbox{\upshape lk}\,}
\newcommand{\Tor}{\mbox{\upshape Tor}}
\newcommand{\astar}{\mbox{\upshape ast}\,}

\newcommand{\lex}{\mbox{\upshape {\tiny lex}}}

\newcommand{\revlex}{\mbox{\upshape {\tiny rl}}}
\newcommand{\Gin}{\mbox{\upshape Gin}}
\newcommand{\field}{{\bf k}}

\newcommand{\MON}{\mbox{\upshape \MON}}
\newcommand{\Init}{\mbox{\upshape in}}

\title{Reverse lexicographic and lexicographic shifting}

\author{Eric Babson}
\author{Isabella Novik}
\author{Rekha Thomas}
\thanks{Research partially supported by NSF grants DMS 0070571 and DMS
0100141} 
\address{Department of Mathematics, University of
  Washington, Seattle, WA 98195-4350, email:
[babson,novik,thomas]@math.washington.edu}   
\date{\today}
\begin{document}

\begin{abstract}
  A short new proof of the fact that all shifted complexes are fixed
  by reverse lexicographic shifting is given.  A notion of
  lexicographic shifting, $\Delta_{\lex}$ --- an operation that
  transforms a monomial ideal of $S=\field[x_i: i\in\N]$ that is
  finitely generated in each degree into a squarefree strongly stable
  ideal --- is defined and studied.  It is proved that (in contrast to
  the reverse lexicographic case) a squarefree strongly stable ideal
  $I\subset S$ is fixed by lexicographic shifting if and only if $I$
  is a universal squarefree lexsegment ideal (abbreviated USLI) of
  $S$.  Moreover, in the case when $I$ is finitely generated and is
  not a USLI, it is verified that all the ideals in the sequence
  $\{\Delta_{\lex}^i(I)\}_{i=0}^{\infty}$ are distinct.  The limit
  ideal
  $\overline{\Delta}(I)=\lim_{i\rightarrow\infty}\Delta_{\lex}^i(I)$
  is well defined and is a USLI that depends only on a certain analog
  of the Hilbert function of $I$.
\end{abstract}

\maketitle

\section{Introduction}

This paper deals with  two problems related to {\em algebraic
shifting} that were raised by Gil Kalai in \cite{K00}.

Algebraic shifting is an algebraic operation introduced
by  Kalai \cite{BK}, \cite{K91} that transforms a simplicial
complex $\Gamma$ into a simpler ({\em shifted}) complex 
$\Delta(\Gamma)$, while preserving
important combinatorial, topological and algebraic invariants such as
face numbers, reduced Betti numbers and extremal algebraic Betti
numbers.  
There are two versions of algebraic shifting ---
exterior and symmetric: the first one amounts to computing 
the (degree) reverse lexicographic generic initial ideal ($\Gin_{\revlex}$) of 
the Stanley-Reisner ideal of $\Gamma$ in the exterior
algebra, while the second one amounts to computing  $\Gin_{\revlex}$
in the symmetric algebra and then applying  a certain ``squarefree''
operation $\Phi$. In this paper we consider only the symmetric version
of algebraic shifting. We refer to  this operation as {\em revlex shifting}
and denote it by $\Delta_{\revlex}$.
%In both cases the resulting ideals 
%are the Stanley-Reisner ideals of shifted simplicial complexes. Those
%complexes are referred to as
% exterior and symmetric (rev-lex)
%algebraic shiftings of $\Gamma$, respectively. 

%Both operations (exterior and symmetric) satisfy the following properties:
%\begin{enumerate}
%\item[(P0)] If $\Gamma$ is a simplicial complex on the vertex set $[n]$,
%then $\Delta(\Gamma)$ is a shifted complex on $[n]$;
%\item[(P1)] $f(\Delta(\Gamma))=f(\Gamma)$;
%\item[(P2)] if $\Gamma'$ is a subcomplex of $\Gamma$, then
%$\Delta(\Gamma')$ is a subcomplex of $\Delta(\Gamma)$;
%\item[(P3)] $\beta_i(\Delta(\Gamma))=\beta_i(\Gamma)$ for
% all $i\geq 0$, where $\beta_i$ stands for the $i$-th reduced Betti number
%of $\Gamma$ computed with coefficients in a certain field $\field$;
%\item[(P4)] $\Delta(\Cone(\Gamma))=
%       \Cone(\Delta(\Gamma))$;
%\item[(P5)] 
%if $\Gamma$ is shifted, then $\Delta(\Gamma)=\Gamma$. In particular,
%$\Delta\Delta\Gamma=\Delta\Gamma$ for any simplicial complex $\Gamma$.
%\end{enumerate}

Clearly $\Delta_{\revlex}(\Gamma)\neq \Gamma$ if $\Gamma$ is not shifted.
Among the many beautiful properties of revlex shifting is the fact that
the converse statement holds as well, namely that
\begin{equation}           \label{P5}
\Delta_{\revlex}(\Gamma)=\Gamma \quad \mbox{if } \Gamma \mbox{ is shifted,}
\end{equation}
and hence that $\Delta_{\revlex}(\Delta_{\revlex}(\Gamma))=
\Delta_{\revlex}(\Gamma)$ for an arbitrary complex $\Gamma$. 
This result was stated in \cite{K91}
and a somewhat hard proof was given in \cite{AHH}.
Eq.~(\ref{P5}) along with the two problems on algebraic shifting
 posed by Gil Kalai \cite[Problems 16 \& 5]{K00} is the starting
point of our paper.

In \cite[Problem 16]{K00} Kalai asks if algebraic shifting
can be axiomatized. In that direction we prove the following result.
(We denote by $[n]$  the set $\{1, 2, \ldots, n\}$, and by $f(\Gamma)$ and 
$\beta_i(\Gamma)$, $i\geq 0$, the $f$-vector and 
the reduced simplicial Betti numbers of $\Gamma$ computed with coefficients
in a field \field, respectively.)
\begin{theorem}     \label{stable}
Let $\Delta$ be an operation  that associates with every $n\geq 0$ 
and every
simplicial complex $\Gamma$ on the vertex set $V=[n]$ a 
shifted simplicial complex $\Delta(\Gamma)$ on the same vertex set.
 Assume further that $\Delta$ satisfies the following properties:
\begin{enumerate}
\item $f(\Delta(\Gamma))=f(\Gamma)$;
\item 
$\Delta(\Gamma\ast\{n+1\})=\Delta(\Gamma)\ast\{n+1\}$;
\item if $\Gamma'\subseteq\Gamma$, then
$\Delta(\Gamma')\subseteq \Delta(\Gamma)$; 
\item 
$\sum_{i=0}^{\dim\Gamma}\beta_i(\Gamma)\leq 
\sum_{i=0}^{\dim\Delta(\Gamma)}\beta_{i}(\Delta(\Gamma))$.
\end{enumerate}
Then for every shifted complex $\Gamma$, $\Delta(\Gamma)=\Gamma$.
\end{theorem}
As a corollary we obtain a new and much simpler proof of Eq.~(\ref{P5}).
(Here $\Gamma\ast\{n+1\}$ is {\em the cone} over $\Gamma$, 
that is, a simplicial
complex on the vertex set $[n+1]$ whose set of faces consists of faces of 
$\Gamma$ together
with $\{F\cup\{n+1\} : F\in \Gamma\}$.)

Problem 5 in \cite{K00} asks whether the
property given by Eq.~(\ref{P5}) holds if
one considers symmetric shiftings with respect to arbitrary
term orders. 
Since in the case of exterior shiftings the answer is positive (as was shown 
by Kalai \cite[Prop.~4.2]{Ka90}),
one may expect to have the same result in the symmetric case as well. 
Here we consider (degree) lexicographic order, 
and denote the corresponding shifting operation by $\Delta_{\lex}$. 
To our surprise we discover that only very few shifted complexes are fixed
by lex shifting.
Our results are summarized in Theorem \ref{Thm2} below. 

Denote by $\N$ the set of all positive integers. We say that an ideal 
$I\subset S=\field[x_i: i\in\N]$ is a 
{\em universal squarefree lexsegment ideal} (abbreviated USLI) if it is 
 finitely generated in each degree and is a squarefree
lexsegment ideal of $S$. 
(Equivalently, an ideal $I$ of $S$ that is finitely generated
in each degree is a USLI if $I\cap S_{[n]}$  is a 
squarefree lexsegment ideal of $S_{[n]}:=\field[x_1, \ldots, x_n]$
for every $n$.)
Thus, for example, the ideal 
$\langle x_1x_2, x_1x_3, x_1x_4x_5x_6x_7\rangle$ is a USLI, while the ideal
$\langle x_1x_2, x_1x_3, x_2x_3\rangle$
is a squarefree lexsegment of $S_{[3]}$ 
but is not a squarefree lexsegment of $S$,
and hence is not a USLI. 
A simplicial complex $\Gamma$ is a {\em USLI complex}
if its Stanley-Reisner ideal, $I_\Gamma$, is a USLI.% (considered as an ideal of $S$).

Recall that for a monomial ideal $J\subset S_{[n]}$
the {\em (bi-graded) Betti numbers} of $J$ are the invariants 
$\beta_{i,j}(J)$ that appear in the minimal free resolution
of $J$ as an $S_{[n]}$-module.
\begin{equation}    \label{alg_betti}
\ldots 
 \bigoplus_j S_{[n]}(-j)^{\beta_{i,j}(J)} \rightarrow \ldots
\rightarrow \bigoplus_j S_{[n]}(-j)^{\beta_{1,j}(J)} \rightarrow
\bigoplus_j S_{[n]}(-j)^{\beta_{0,j}(J)} \rightarrow J \rightarrow 0
\end{equation}
Here $S_{[n]}(-j)$ denotes $S_{[n]}$ with grading shifted by $j$.
Following \cite{Eis2}, we define the {\em $B$-sequence} of $J$,
$B(J):=\{B_j(J) : j\geq 1\}$, where $B_j(J):=\sum_{i=0}^j (-1)^i
\beta_{i,j}(J)$. (The $B$-sequence of an ideal contains the same
information as its Hilbert series --- see Section
\ref{infinite_section} for more details as well as for the definition
of the $B$-sequence for a monomial ideal of $S$ that is finitely
generated in each degree.)

\begin{theorem} \qquad       \label{Thm2}
\begin{enumerate} 
\item A (finite) shifted simplicial complex $\Gamma$ satisfies
$\Delta_{\lex}(\Gamma)=\Gamma$ if and only if $\Gamma$ is a USLI complex.
 Moreover, if $\Gamma$
is not a USLI complex, then  all the complexes
in the sequence $\{\Delta_{\lex}^i(\Gamma)\}_{i=0}^{\infty}$ are distinct.
(Here $\Delta_{\lex}^i(\Gamma)$ denotes the complex obtained from 
$\Gamma$ by $i$
consecutive applications of $\Delta_{\lex}$.)

\item The {\em limit} complex
$\overline{\Delta}_{\lex}(\Gamma):=
\lim_{i\rightarrow\infty}\Delta_{\lex}^i(\Gamma)$
 is well defined and is a (usually infinite) USLI complex.
 Moreover, $\overline{\Delta}_{\lex}(\Gamma)$ is 
the unique USLI complex whose Stanley-Reisner ideal has the same $B$-sequence 
as $I_\Gamma$.
\end{enumerate}
\end{theorem}

The last part of the theorem implies that
if two simplicial complexes $\Gamma_1$ and $\Gamma_2$
that have the same $h$-vector 
(up to possibly several zeros appended at the end), then
$\overline{\Delta}_{\lex}(\Gamma_1)=\overline{\Delta}_{\lex}(\Gamma_2)$.
Thus, in contrast to revlex shifting,
 the operation $\overline{\Delta}_{\lex}$
forgets all the information that $\Gamma$ carries (including the dimension
of $\Gamma$) except its $h$-numbers.

Our theorems establish for simplicial complexes, results similar in
spirit to those in commutative algebra due to Bigatti-Conca-Robbiano
\cite{BCR} and Pardue \cite{Pardue}.
%  on the behavior of (degree)
% reverse lexicographic and (pure) lexicographic generic initial ideals
% with respect to distractions: 
Theorem~4.3 in \cite{BCR} asserts that if $I$ is a strongly stable
ideal in $S_{[n]}$ and $\mathcal{L}$ is a distraction matrix, then
$\Gin_{\revlex}(D_{\mathcal L}(I))=I$, while Proposition~30 in
\cite{Pardue} asserts that sufficiently (but finitely) many
applications of the operation $\Gin_{\lex} \circ D_{\mathcal L}$ to a
monomial ideal $I\subset S_{[n]}$ results in the unique lexsegment
ideal of $S_{[n]}$ having the same Hilbert function as $I$.

The structure of the paper is as follows. Section \ref{axiom_section}
is devoted to the proof of Theorem \ref{stable}. 
In Section \ref{revlex_section}
after recalling basic facts and definitions related to generic initial ideals
and revlex shifting we provide a short new proof of Eq.~(\ref{P5}). 
In Section \ref{USLI_section} we introduce and study the class of
universal squarefree lexsegment ideals (USLIs) 
and the class of {\em almost USLIs} ---
the notions that play a crucial role in the proof of Theorem \ref{Thm2}.
 Finally in Section \ref{infinite_section} we prove Theorem \ref{Thm2}.
We close with a brief discussion of arbitrary term orders.

\section{Axiomatizing Algebraic Shifting}    \label{axiom_section}
This section is devoted to the proof of Theorem \ref{stable}.
 We start by reviewing several
notions pertaining to simplicial complexes.

Denote the collection of all subsets of 
$[n]:=\{1, 2, \ldots, n\}$ by $2^{[n]}$.
Recall that a simplicial complex $\Gamma$ on the vertex set $[n]$ 
is a collection $\Gamma\subseteq 2^{[n]}$ that  is
closed under inclusion. 
(We do not require that every singleton $\{i\}\subseteq[n]$
is  an element of $\Gamma$.)
The elements of $\Gamma$ are called faces
and the maximal faces (under inclusion) are called {\em facets}.
 $F\in\Gamma$ is an
{\em $i$-dimensional face} (or an $i$-face) if
$|F|=i+1$. The {\em dimension} of $\Gamma$,
$\dim \Gamma$, is the maximal dimension of its faces.
 The number of $i$-dimensional faces of $\Gamma$
is denoted by $f_i(\Gamma)$, and the sequence 
$f(\Gamma):=(f_{-1}(\Gamma), f_0(\Gamma), f_1(\Gamma), 
\ldots f_{\dim\Gamma}(\Gamma))$ is called the
$f$-vector of $\Gamma$. 
Another set of invariants associated
with $\Gamma$ is the set of its  reduced
Betti numbers 
$\beta_i(\Gamma):=\dim_\field\widetilde{H}_i(\Gamma; \field)$,
where $\widetilde{H}_i(\Gamma; \field)$ is the $i$-th reduced
simplicial homology of $\Gamma$ with coefficients in a field $\field$.

A simplicial complex $\Gamma$ on the vertex set 
$[n]$ is called {\em shifted} if for every
$F\in\Gamma$, $i\in F$, and $i<j\leq n$, the set 
$(F\setminus\{i\})\cup\{j\}$ is a face of $\Gamma$ as well.
The Betti numbers of a shifted complex $\Gamma$ are combinatorial
invariants and can be  computed via the 
following well-known formula \cite[Thm.~4.3]{BK}:
\begin{lemma}   \label{betti}
If $\Gamma$ is a shifted complex on the vertex set $[n]$, then 
$$\beta_i(\Gamma)=|\{F\in\max(\Gamma) \,:\, |F|=i+1,\, n\notin \Gamma\}|,$$
where $\max(\Gamma)$ denotes the set of 
facets of $\Gamma$.
\end{lemma}

For a simplicial complex $\Gamma$ and a vertex $v$ of $\Gamma$ define
the {\em antistar of $v$ in $\Gamma$} as
$\astar_{\Gamma} (v) = \{F\in\Gamma: v \notin F\}$. Define also 
the {\em link of $v$ in $\Gamma$} by $\lk_\Gamma (v):=
\{F\in\astar_{\Gamma}(v) :  F\cup\{v\}\in\Gamma\}$.
Note that if $\Gamma$ is a shifted complex on the vertex set $[n]$,
then $\lk_\Gamma (n)$ and $\astar_\Gamma (n)$ are 
shifted complexes on $[n-1]$. 

If $\Gamma$ is a simplicial complex on $V$ and $u\not\in V$,
then the {\em cone} over $\Gamma$ with apex $u$ is a simplicial complex,
denoted  $\Gamma\ast \{u\}$, on the vertex set $V\cup\{u\}$
whose faces are all sets of the form   
$F\cup A$, where $F\in\Gamma$ and $A\subseteq \{u\}$.
Thus 
for any vertex $v$ of $\Gamma$, 
$\Gamma=\lk_\Gamma(v)\ast\{v\} \cup \astar_\Gamma(v)$ and 
$\lk_\Gamma(v)\ast\{v\} \subseteq \Gamma \subseteq 
             \astar_\Gamma(v)\ast\{v\}$.

Now we are ready to verify Theorem \ref{stable} asserting that
if $\Delta$ is an operation  that associates with every $n\geq 0$ 
and every
simplicial complex $\Gamma$ on the vertex set $V=[n]$ a 
shifted simplicial complex $\Delta(\Gamma)$ on the same vertex set,
and if
$\Delta$ satisfies the following properties:
\begin{enumerate}
\item $f(\Delta(\Gamma))=f(\Gamma)$;
\item 
$\Delta(\Gamma\ast\{n+1\})=\Delta(\Gamma)\ast\{n+1\}$;
\item if $\Gamma'\subseteq\Gamma$, then
$\Delta(\Gamma')\subseteq \Delta(\Gamma)$; 
\item 
$\sum_{i=0}^{\dim\Gamma}\beta_i(\Gamma)\leq 
\sum_{i=0}^{\dim\Delta(\Gamma)}\beta_{i}(\Delta(\Gamma))$,
\end{enumerate}
then for every shifted complex $\Gamma$, $\Delta(\Gamma)=\Gamma$.

\smallskip\noindent {\it Proof of Theorem \ref{stable}: \ } 
Fix a shifted complex $\Gamma$ on $n$ vertices.
If $n=0$ or $n=1$ then $\Delta(\Gamma)=\Gamma$ by property (1). 
We proceed by induction on~$n$. Since
 the link and the antistar of the vertex $n$
in $\Gamma$, $\Gamma'=\lk_\Gamma (n)$ and 
$\Gamma''=\astar_\Gamma (n)$, 
respectively,
are shifted complexes on the vertex set $[n-1]$
and since 
$\Gamma'\ast \{n\}\subseteq \Gamma \subseteq \Gamma''\ast\{n\}$,
the induction hypothesis together with
properties (2) and (3) yield
$$\Gamma'\ast \{n\}\subseteq \Delta(\Gamma) 
\subseteq \Gamma''\ast\{n\}.$$
Therefore,
\begin{eqnarray*}
A &:=& \{F\in\max(\Delta(\Gamma)) \,:\, n\notin F\} 
=\{F\in\max(\Delta(\Gamma)) \,:\,  F\in \Gamma''\} \\
&\subseteq& \Gamma''\setminus\Gamma' \stackrel{(\star)}{=}
\{F\in\max(\Gamma) \,:\, n\notin F\}=:B,
\end{eqnarray*}
where 
$(\star)$ follows from the shiftedness of $\Gamma$:
$$F\in\max(\Gamma)\cap 2^{[n-1]}  \Longleftrightarrow  
F\in\Gamma \mbox{ but } F\cup\{n\} \notin \Gamma \Longleftrightarrow 
F\in \astar_\Gamma(n)\setminus \lk_\Gamma(n)=\Gamma''\setminus \Gamma'.$$ 
On the other hand, Lemma \ref{betti} and property (4) imply that
$$|A|=\sum_{i=0}^{\dim\Delta(\Gamma)}\beta_i(\Delta(\Gamma))\geq 
\sum_{i=0}^{\dim\Gamma} \beta_i(\Gamma)=|B|,$$
and thus that $A=B$. Hence
$\Delta(\Gamma)\supseteq A=\Gamma''\setminus \Gamma'$, and
we infer that 
$$\Delta(\Gamma)\supseteq (\Gamma'\ast\{n\})\cup( \Gamma''\setminus \Gamma')=
\Gamma.$$
 Since $f(\Gamma)=f(\Delta(\Gamma))$ by property (1), it follows
that $\Delta(\Gamma)=\Gamma$. \endproof

\section{Generic Initial Ideals and revlex shifting}          
                        \label{revlex_section}
In this section we review basic facts and definitions related to
generic initial ideals and revlex shifting.  We also provide a new short proof
of Eq.~(\ref{P5}) asserting that
 $\Delta_{\revlex}(\Gamma)=\Gamma$ for a shifted 
$\Gamma$.
Let $S_{[n]}=\field[x_1, \ldots, x_n]$ 
be the ring of polynomials in $n$  variables
 over  a field $\field$ of characteristic zero, and let $\Gamma$ be a
 simplicial complex on the vertex set $[n]$. We recall that
the {\em Stanley-Reisner ideal} of  $\Gamma$ 
\cite{St} is the squarefree monomial  ideal $I_\Gamma \subset S_{[n]}$ 
whose generators correspond to nonfaces of $\Gamma$:
$$ I_{\Gamma} := \langle \prod_{j=1}^k x_{i_j} \in S_{[n]}\,:\,
\{i_1<i_2<\ldots<i_k\}\notin\Gamma\rangle.
$$
The Stanley-Reisner ideal of a shifted complex is called a
{\em squarefree strongly stable ideal}. 
(Equivalently, a squarefree monomial ideal  
$I$ is squarefree strongly stable, if for every minimal generator $m$ of 
$I$ and for
every $1\leq i<j$ such that $x_j|m$ but $x_i\not|m$, the monomial $mx_i/x_j$
lies in $I$.)

Let $\succ$ be a term order on $S_{[n]}$ that refines the partial
order by degree where lower degree monomials are more expensive than
higher degree monomials, and satisfies $x_1\succ\ldots\succ x_n$. Let
$I\subset S_{[n]}$ be a homogeneous ideal such as the Stanley-Reisner
ideal of $\Gamma$.  Consider a generic $n\times n$ matrix $g$.  Then
$g$ acts on the set of linear forms of $S_{[n]}$ by
$gx_j=\sum_{i=1}^{n}g_{ij}x_i$ and this action can be extended
uniquely to a ring automorphism on $S_{[n]}$ that we also denote by
$g$.  Following \cite[Thm.~15.18]{Eis} define the {\em generic initial
  ideal of} $I$ with respect to $\succ$ as
$$\Gin_\succ(I):=\Init_\succ(gI),$$
where $\Init_\succ(gI)$ is the
initial ideal of $gI$ with respect to $\succ$ in the sense of
Gr\"obner basis theory. The same theorem in \cite{Eis} asserts that we
can choose $g$ to be upper triangular and hence we assume from now on
that $gx_j=\sum_{i=1}^jg_{ij}x_i$.

We briefly outline how to compute $\Gin_\succ(I)$ (for a detailed description
the reader is referred to \cite[Thm.~15.18]{Eis}).
\begin{definition}
An exterior monomial in  $\bigwedge^l(S_{[n]})_d$ is an element of the form 
$m_1\wedge\ldots\wedge m_l$ where each $m_i$ is a monomial of $S_{[n]}$ 
of degree $d$ and $m_1\succ \ldots \succ m_l$. The 
extension of $\succ$ to an order on monomials
of $\bigwedge^l(S_{[n]})_d$ is the order in which 
$m_1\wedge\ldots\wedge m_l \succ n_1\wedge\ldots\wedge n_l$ if for some $s$
we have that $m_s\succ n_s$ and $m_i=n_i$ for $i<s$. For a non-zero element 
$f$ of $\bigwedge^l(S_{[n]})_d$, define the {\em initial term}  of $f$, 
$\Init_\succ(f)$, to be
the $\succ$-largest monomial appearing in $f$ with nonzero coefficient
when $f$ is written as a linear combination of (distinct) monomials.
\end{definition}

Consider a generic $n\times n$ upper-triangular matrix $g$ 
and its action on $S_{[n]}$.
Let $I_d$ be the $d$-th homogeneous component of a homogeneous ideal $I$, 
and let  $f_1, \ldots,  f_t$ be a basis of $I_d$.
Then $g(f_1)\wedge \ldots \wedge g(f_t)\in\wedge^t(S_{[n]})_d$.
Denote by  $M_d=m_1\wedge \ldots \wedge m_t$ the monomial
$\Init_\succ (g(f_1)\wedge \ldots \wedge g(f_t))$ 
%--- the $\succ$-largest
%monomial appearing in $g(f_1)\wedge \ldots \wedge g(f_t)$, 
and by $V_d$  the subspace of $(S_{[n]})_d$
spanned by $m_1, \ldots, m_t$.

\begin{proposition}  \label{gin_constr}
 $\Gin_\succ(I)=\bigoplus V_d$.
\end{proposition}

Several basic properties of Gins are summarized in the following lemma.
\begin{lemma}  \label{cone}
Let $I\subset S_{[n]}$  be a homogeneous ideal. Then 
\begin{enumerate}
\item $\Gin_\succ(I)$ is a strongly stable monomial ideal (that is, if 
$m\in \Gin_\succ(I)$, $x_j|m$ and $1\leq i <j$, then 
$x_im/x_j\in \Gin_\succ(I)$ as well). 
\item $\Gin_\succ(I)$ and $I$ have the same Hilbert function 
(that is, $\dim_\field(\Gin_\succ(I)_d)=\dim_\field(I_d)$ for all $d$).
\item If $J\subseteq I$ is a homogeneous ideal of $S_{[n]}$, then 
$\Gin_\succ(J)\subseteq \Gin_\succ(I)$.
\item Let
$\succ'$ be an extension of  $\succ$ to a term order on $S_{[n+1]}$ satisfying
$x_n\succ' x_{n+1}$. Then 
$\Gin_{\succ'}(IS_{[n+1]})=(\Gin_\succ I) S_{[n+1]}.$
In particular, for a simplicial complex $\Gamma$ on $[n]$,
$\Gin_{\succ'}(I_{\Gamma\ast\{n+1\}})=(\Gin_\succ I_\Gamma) S_{[n+1]}.$
\end{enumerate}
\end{lemma}
\proof
Part (1)  is \cite[Thm.~15.18 and Thm.~15.23]{Eis}. 
Part (2) follows from \cite[Thm 15.3]{Eis}.
Part (3) is obvious from the definitions.
To  prove part (4),  consider a generic upper-triangular 
$(n+1)\times(n+1)$ matrix  $\widetilde{g}$ and 
its left-upper $n\times n$ submatrix $g$. 
Then $g$ acts on $S_{[n]}$, $\widetilde{g}$
acts on $S_{[n+1]}$, and 
$\widetilde{g}x_i=gx_i$ for all $1\leq i\leq n$. 
Therefore for every (homogeneous)  element 
$h$ of $I\subset S_{[n]}\subset S_{[n+1]}$,
$\widetilde{g}h=gh$. Thus for $h\in I$,
 $\Init_{\succ'}(\widetilde{g}h)=\Init_{\succ}(gh)$, implying that
$\Gin_\succ I \subseteq \Gin_{\succ'}(IS_{[n+1]})$, and hence that
$(\Gin_\succ I)S_{[n+1]}\subseteq \Gin_{\succ'}(IS_{[n+1]})$. 
The lemma follows, since both  the ideals $(\Gin_\succ I)S_{[n+1]}$ and 
$\Gin_{\succ'}(IS_{[n+1]})$
have the same Hilbert function.  \endproof

In the later sections we compare Gins of the same ideal $I$
computed with respect to different term orders. 
For that  we need the following definition:

\begin{definition}
  Let $I_1\neq I_2$ be two monomial ideals of $S_{[n]}$ and let $\succ$ be
  a term order.  We say that $I_1\succ I_2$ if the largest monomial in
  the symmetric difference of $I_1$ and $I_2$ is in $I_1$.
  Equivalently, $I_1\succ I_2$ if the largest monomial in the
  symmetric difference of $G(I_1)$ and $G(I_2)$ is in $G(I_1)$, where
  $G(I_1)$ and $G(I_2)$ are the sets of minimal generators of $I_1$
  and $I_2$ respectively.
\end{definition}

One immediate observation is
\begin{lemma}  \label{gin>gin}
Let $\sigma$ and $\tau$ be two term orders on $S_{[n]}$.
Then $\Gin_{\sigma}(I)\geq_\sigma \Gin_{\tau}(I)$
for any homogeneous ideal $I\subset S_{[n]}$.
\end{lemma} 
\proof
Let $f_1, \ldots, f_t$ be a basis of $I_d$, and let 
$g$ be a generic $n \times n$
upper-triangular matrix. 
Since $M'_d:=\Init_{>_{\tau}}(g(f_1)\wedge\ldots\wedge g(f_t))$ appears
in $g(f_1)\wedge\ldots\wedge g(f_t)$ with a non-zero coefficient, it follows
that 
$M_d:=\Init_{\sigma}(g(f_1)\wedge\ldots\wedge g(f_t)) \geq_{\sigma} M'_d$
(for every $d\geq 0$).
Proposition \ref{gin_constr} implies the lemma.
\endproof

We remark that a stronger version of Lemma \ref{gin>gin} was proved in 
\cite[Cor.~1.6]{Conca}.

Another ingredient needed for defining revlex shifting
is the notion of the squarefree operation. This is a bijection $\Phi$
between the set of all monomials in $\{x_i : i\in \N\}$ and the set of all
squarefree monomials in $\{x_i : i\in \N\}$, defined by
$$\Phi(\prod_{j=1}^k x_{i_j})=\prod_{j=1}^k x_{i_j+j-1}, \mbox{ where }
   i_1\leq i_2\leq \ldots\leq i_k.
$$
 Note that for a monomial  $m\in S_{[n]}$, 
$\Phi(m)$ may not belong to  $S_{[n]}$. However the graded 
reverse lexicographic order 
 has the following remarkable property 
\cite[Lemma 6.3(ii)]{K91}, \cite[Lemma 1.1]{AHH}: if 
$m$ is a minimal generator of $\Gin_{\revlex} I_\Gamma$ 
(where $\Gamma$ is a simplicial 
complex on $[n]$), then $\Phi(m)$ is an element of $S_{[n]}$.
This leads to the following definition (due to Kalai):
\begin{definition}   \label{equiv}
Let $\Gamma$ be a simplicial complex on the vertex set $[n]$.
The reverse lexicographic shifting of $\Gamma$,
$\Delta_{\revlex}(\Gamma)$, is a simplicial  complex on $[n]$ whose 
Stanley-Reisner ideal is given by
$$I_{\Delta_{\revlex}(\Gamma)}=\langle \Phi(m) \; : \; m\in
G(\Gin_{\revlex} I_\Gamma) \rangle,$$
where $G(I)$ denotes the set of
the minimal generators of a monomial ideal $I$.
\end{definition}

We now provide a new and simple proof of Eq.~(\ref{P5})
 (due originally to Aramova, Herzog, and Hibi \cite{AHH}).
\begin{theorem}   \label{AHH}
The revlex shifting $\Delta_{\revlex}$ satisfies all the conditions of 
Theorem \ref{stable}. Thus $\Delta_{\revlex}(\Gamma)=\Gamma$ for every
 shifted complex $\Gamma$.
\end{theorem} 
\proof It is well-known that (symmetric) revlex shifting satisfies all
the conditions of Theorem \ref{stable}, except possibly for property
(2) whose proof appears to be missing in the literature (for the
exterior version of algebraic shifting it was recently verified by
Nevo \cite{Nevo}): the fact that $\Delta(\Gamma)$ is a shifted
simplicial complex follows from Lemma~\ref{cone}(1); property (1) is
\cite[Lemma 6.3(i)]{K91}; property (3) is a consequence of 
Lemma~\ref{cone}(3); property (4) follows from \cite[Cor.~8.25]{H} asserting
that $\beta_i(\Gamma)=\beta_i(\Delta(\Gamma))$ for all $i$. To prove
property (2) it suffices to check that $\Delta(\Gamma)$ and
$\Delta(\Gamma\star \{n+1\})$ have the same set of minimal nonfaces
(equivalently, $I_{\Delta(\Gamma)}\subset S_{[n]}$ and
$I_{\Delta(\Gamma\star \{n+1\})}\subset S_{[n+1]}$ have the same set
of minimal generators).  This follows from Definition \ref{equiv} and
Lemma \ref{cone}(4). \endproof

%\begin{remark} \quad
\paragraph{\bf Remarks} \quad

(1)
  We note that to verify the inequality $\sum \beta_i(\Gamma)\leq \sum
  \beta_i(\Delta_{\revlex}(\Gamma))$ one does not need to use the fact
  that $\beta_i(\Gamma)=\beta_i(\Delta_{\revlex}(\Gamma))$ for all
  $i$, which is a consequence of the deep result due to
  Bayer--Charalambous--Popescu \cite{BCP} and Aramova--Herzog~\cite{AH} 
that revlex shifting preserves extremal (algebraic) Betti
  numbers.  Instead one can use the standard flatness argument (see
  \cite[Thm.~3.1]{H}) to show that $\beta_{i,j}(I_\Gamma) \leq
  \beta_{i,j}(\Gin_{\revlex}(I_\Gamma)) = \beta_{i,j}(I_{\Delta(\Gamma)})$ 
for all $i$, $j$, where the equality comes from
  the fact that $\Phi$ applied to (minimal generators of)
 a strongly stable ideal
  $\Gin_{\revlex}(I_\Gamma)$ preserves algebraic Betti numbers (see
  \cite[Lemma 2.2]{AHH}). The Hochster formula \cite{Hoc} then asserts
  that the reduced Betti numbers of a simplicial complex are equal to
  certain algebraic graded Betti numbers of its Stanley-Reisner ideal.

(2)
 In algebraic terms, the statement of
Theorem \ref{AHH} translates to the fact that if
$I\subset S_{[n]}$ is a squarefree strongly stable ideal, then
$\Phi(\Gin_{\revlex}(I))=I$, where  
$\Phi(\Gin_{\revlex}(I)):=\langle \Phi(m): m \in G(\Gin_{\revlex}(I)) \rangle
$.
Hence $\Gin_{\revlex}(I)=\langle\Phi^{-1}(\mu) : \mu\in G(I)\rangle$, that is,
 computing the revlex Gin of a squarefree strongly stable ideal $I$
simply amounts to applying $\Phi^{-1}$ to the  minimal generators
of $I$.

(3) Our proof (as well as the original proof  in \cite{AHH})
of the equation $\Phi(\Gin_{\revlex}(I))=I$ for a squarefree strongly
stable ideal $I$ works only over a field $\field$ of characteristic zero.
We however do not know of any counterexamples in the case of a field
of positive characteristic.
%\end{enumerate}
%\end{remark}

\section{Combinatorics of USLIs, almost USLIs, and lex Gins}   
                 \label{USLI_section}
In this section we introduce and study 
the class of universal squarefree lexsegment ideals (USLIs)
and the class of almost USLIs. These notions turn out to be crucial 
in the proof of Theorem \ref{Thm2}.
To allow for infinitely generated ideals 
(as we  need in the following section) we 
consider the system of rings 
$S_{[n]}$, $n\in\N$,
endowed with natural embeddings $S_{[n]}\subseteq S_{[m]}$ for $m\geq n$, and
provide definitions suitable
for the direct limit 
ring $S=\lim_{n\rightarrow\infty}S_{[n]}=\field[x_i : i\in\N]$.

Recall that a squarefree monomial ideal 
$I\subset S$ ($I\subset S_{[n]}$, respectively) is a 
{\em squarefree lexsegment ideal} of $S$ ($S_{[n]}$, respectively)
 if for every
monomial $m\in I$ and every
squarefree monomial $m'\in S$ ($m'\in S_{[n]}$, respectively)
such that $\deg(m')=\deg(m)$ and $m'>_{\lex} m$, 
$m'$ is an element of $I$ as well.
\begin{definition}   \label{USLI}
An ideal $L$ of $S$ (or of $S_{[n]}$)
is
a {\em universal squarefree lexsegment ideal} (abbreviated USLI)
if it is finitely generated in each degree and $LS$ is
a squarefree lexsegment ideal of $S$. 
Equivalently, an ideal 
$L=L(k_\bullet)$ (here $k_\bullet=\{k_i\}_{i\in\N}$ is a 
sequence of nonnegative integers) is a 
USLI with $k_i$  minimal generators
of degree $i$ (for $i\in\N$) if and only if
 the set of minimal generators of $L$, 
$G(L)$, is given by
$$G(L)=
\bigcup_{r=1}^{\infty}\left\{
(\prod_{j=1}^{r-1} x_{R_j})\cdot x_l \;:\; R_{r-1}+1 
\leq l \leq R_r-1\right\},  \mbox{ where } R_j=j+\sum_{i=1}^j k_i. 
  %         0\leq j \leq \infty. 
$$
\end{definition}

The easiest way to verify the description of the set
$G(L)=\{m_1 >_{\lex} m_2  >_{\lex} \cdots \ >_{\lex} m_s  >_{\lex} \cdots \}$ 
of a USLI $L$ is by induction on $s$. Indeed, if $m_1, \cdots, m_s$ satisfy
the above description and
$m_s=(\prod_{j=1}^{r-1} x_{R_j})\cdot x_l$, 
then there are two possibilities for
$m_{s+1}$: either $\deg(m_{s+1})=\deg(m_s)=r$ (equivalently, $l<R_r-1$)
or $\deg(m_{s+1})=r'>r$ (equivalently, $l=R_r-1$ and $k_i=0$ for all $r<i<r'$).
In the former case, since $m_{s}>_{\lex} m_{s+1}$ and since $m_s$ is the 
immediate lex-predecessor of
$m':=(\prod_{j=1}^{r-1} x_{R_j})\cdot x_{l+1}$,
it follows that
$m'\geq_{\lex} m_{s+1}\in L$ 
which together
with  $L$ being a USLI implies that 
$m'\in L$. Since $m'$ is not divisible by any of $m_1, \cdots, m_s$,
this yields
$m_{s+1}=m'$. The treatment
of the latter case is similar: just observe that every squarefree monomial
of degree $r'$ that is lex-smaller than
 $m':=(\prod_{j=1}^{r-1} x_{R_j})\cdot (\prod_{j=1}^{r'-r+1} x_{l+j})=
 (\prod_{j=1}^{r'-1} x_{R_j})\cdot x_{R_{r'-1}+1}$ is divisible by at least
one of
$m_1, \ldots, m_s$ and hence is in $L-G(L)$, while $m'$ is not divisible
by any of $m_1, \cdots, m_s$.

\begin{example} \quad
\begin{enumerate} 
\item
  The ideal $\langle x_1x_2, x_1x_3,
  x_2x_3\rangle$ (the Stanley-Reisner ideal of three isolated points)
  is a lexsegment in $S_{[3]}$, but is not a lexsegment in $S$, and
  hence is not a USLI.  

\item
  The ideal $I = \langle x_1x_2, x_1x_3,
  x_1x_4x_5x_6x_7\rangle$ is the USLI with $k_1 = 0, k_2 = 2, k_3 =
  k_4 = 0, k_5 = 1$ and $k_i = 0$ for all $i > 5$. In this example,
  check that $R_1 = 1, R_2 = 4, R_3 = 5, R_4 = 6$ and $R_5 = 8$. 
\end{enumerate}
\end{example}

Note that every USLI is a squarefree strongly stable ideal, and hence
is the Stanley-Reisner ideal of a shifted (possibly infinite)
simplicial complex (we refer to such complex as a {\em USLI complex}).
All complexes considered in this section are assumed to be finite.

The following lemma describes  certain combinatorial properties of 
USLI complexes. This lemma together with
 Lemmas \ref{main_lemma} and  \ref{Pardue_lemma}
below provides a key step in the proof of 
Theorem \ref{Thm2}.
\begin{lemma}   \label{comb_USLI}
Let $\Gamma$ be a USLI complex on the vertex set $[n]$
with $I_\Gamma=L(k_\bullet)$. 
\begin{enumerate}
\item If $I_\Gamma\neq 0$ and $k_d$ is the last nonzero entry
in the sequence $k_\bullet$, then $\Gamma$ has exactly $d$ facets. 
They are given by
$$F_i=\left\{\begin{array}{ll}
    \{R_j : 1\leq j \leq i-1\}\cup [R_i+1, n] &  \mbox{ if  $1\leq
    i\leq d-1$,}\\ 
    \{R_1, \ldots, R_{d-1}\}\cup [R_d, n] & \mbox{ if $i=d$.}
\end{array}\right.
$$
\item If $\Gamma'$ is a shifted complex on $[n]$ such that 
$f(\Gamma)=f(\Gamma')$, then $\Gamma=\Gamma'$. (In other words
every USLI complex is the only shifted complex in its $f$-class).
\end{enumerate}
\end{lemma} 
\proof We verify part (1) by induction on $n+d+\sum k_i$. The assertion
clearly holds if $d=1$ or if $\sum k_i=1$. For instance, if $d=1$ and $k_1=n$
(equivalently, $R_1=n+1$), then $F_1=[n+1, n]=\emptyset$ is the only 
facet of $\Gamma$.

Note that $R_d$
is the index of the first variable that does
not divide any of the minimal generators of $I_\Gamma$.
Thus if $R_d\leq n$, then $\Gamma=\lk_\Gamma(n)\star\{n\}$, and we are done
by applying induction hypothesis to the USLI complex $\lk_\Gamma(n)$.
So assume that $R_d=n+1$.
Then $\lk_\Gamma(n)$ and $\astar_\Gamma(n)$
are easily seen to be the USLI complexes on the vertex set $[n-1]$
whose Stanley-Reisner ideals are given
by $L_1=L(k_1, \ldots, k_{d-2}, k_{d-1}+1)$ and 
$L_2=L(k_1, \ldots, k_{d-1}, k_d-1)$,
respectively. 
Hence by induction hypothesis the complex
$\lk_\Gamma(n)\star\{n\}$ has exactly $d-1$ facets, namely
the sets $F_1, \ldots, F_{d-1}$ from the list above.
Now if $k_d>1$, then  by induction hypothesis
the facets of $\astar_\Gamma(n)$ are the sets $F_1-\{n\}, \ldots,
F_{d-1}-\{n\}, F_d$. Since
$\Gamma= (\lk_\Gamma(n)\star\{n\}) \cup \astar_\Gamma(n)$,
it follows that  $\max(\Gamma)=\{F_1, \ldots, F_d\}$.
Similarly, if $k_d=1$ and $k_j$ is the last nonzero
entry in the sequence $(k_1, \ldots, k_{d-1})$, then 
 the facets of $\astar_\Gamma(n)$ are the sets $F_1-\{n\}, \ldots,
F_{j-1}-\{n\}, F_d$, and the result follows in this case as well.

To prove part (2) we induct on $n$. The assertion is obvious for $n=1$.
For $n>1$ we consider two cases.

{\bf Case 1:} $R_d\leq n$. In this case 
$\Gamma=\lk_\Gamma(n)\star\{n\}$, so $\beta_i(\Gamma)=0$ for all $i$.
Since among all squarefree strongly stable
ideals with the same Hilbert function
the squarefree lexsegment ideal has the largest 
algebraic Betti numbers \cite[Thm.~4.4]{AHHlex}, 
and since by Hochster's formula \cite{Hoc}, 
$\beta_{n-i-1}(\Lambda)=\beta_{i-1,n}(I_\Lambda)$ 
for any simplicial complex $\Lambda$ on the vertex set $[n]$,
 it follows that $\beta_i(\Gamma')\leq\beta_i(\Gamma)=0$,
and so $\beta_i(\Gamma')=0$ for all $i$. Since $\Gamma'$ is shifted,
Lemma \ref{betti} implies that all facets of $\Gamma'$ contain $n$.
Thus $\Gamma'=\lk_{\Gamma'}(n)\star\{n\}$, and the assertion follows
from induction hypothesis applied to 
$\lk_\Gamma(n)$ and $\lk_{\Gamma'}(n)$.

{\bf Case 2:} $R_d=n+1$. In this case all facets of $\Gamma$ but $F_d$ 
contain vertex $n$ (this follows from part (1) of the Lemma), and we infer
from Lemma \ref{betti}  that
$$\beta_i(\Gamma)=\left\{\begin{array}{ll} 0, \mbox{ if $i\neq d-2$} \\
                                          1, \mbox{ if $i= d-2$.} \\
                        \end{array}
                 \right.
$$
Recall the Euler-Poincar\'e formula asserting that for any simplicial
complex $\Lambda$, 
$$\sum_{j\geq -1}(-1)^j f_j(\Lambda) 
= \sum_{j\geq -1}(-1)^j \beta_j(\Lambda)
=:\widetilde{\chi}(\Lambda).$$
Therefore, $\widetilde{\chi}(\Gamma')=\sum_{j\geq -1}(-1)^j f_j(\Gamma')=
\sum_{j\geq -1}(-1)^j f_j(\Gamma)=\widetilde{\chi}(\Gamma)=(-1)^{d-2}$, and
hence not all Betti numbers of $\Gamma'$  vanish. The
same reasoning as in Case 1 then shows that  
$\beta_i(\Gamma')=\beta_i(\Gamma)$ for all $i$. Applying
Lemma \ref{betti} once again, we obtain that 
$\Gamma'=(\lk_{\Gamma'}(n)\star\{n\})\cup\{ F'\}$, where $|F'|=d-1$
and $F'$
is the only facet of $\Gamma'$ that does not contain $n$. Thus 
$f(\lk_{\Gamma}(n))= f(\lk_{\Gamma'}(n))$ and
 $f(\astar_{\Gamma}(n))= f(\astar_{\Gamma'}(n))$, and so
$\lk_{\Gamma}(n)=\lk_{\Gamma'}(n)$ and 
$\astar_{\Gamma}(n)=\astar_{\Gamma'}(n)$ (by induction hypothesis),
yielding that $\Gamma=\Gamma'$.
\endproof

We now turn to the class of {\em almost USLIs}. 
(Recall our convention that lower degree monomials are 
lex-larger than  higher degree monomials.)

\begin{definition}
Let $I\subset S$ (or $I\subset S_{[n]}$)
be a squarefree strongly stable monomial  ideal with 
$G(I)=\{m_1>_{\lex}  \ldots >_{\lex} m_l>_{\lex} m_{l+1} \}$. 
We say that $I$ is {\em an almost USLI}
if  $I$ is not a USLI, but $L=\langle m_1, \ldots, m_l\rangle$ is a USLI.
We say that a simplicial complex $\Gamma$ is {\em an almost USLI complex} 
if $I_\Gamma$ is an almost USLI.
\end{definition}

As we will see in the next section (see also Lemma \ref{Pardue_lemma} below),
 what makes almost USLI complexes noninvariant under
lex shifting is the following combinatorial property. (We recall that the 
{\em regularity} of a finitely generated stable monomial ideal $I$, $\reg(I)$,
is the maximal degree of its minimal generators.)

\begin{lemma}  \label{main_lemma}
Let $\Gamma$ be an almost USLI complex.
Then $|\max(\Gamma)|>\reg(I_\Gamma)$.
\end{lemma}

\proof 
Assume $\Gamma$ is a simplicial complex on $[n]$
with $G(I_\Gamma)=\{m_1>_{\lex}\ldots>_{\lex}m_l>_{\lex}m_{l+1}\}$.
We have to show that $|\max(\Gamma)|>\deg(m_{l+1})=:d$.
We verify this by induction on $d$. To simplify the notation
 assume without loss of generality  that every singleton 
$\{i\}\subset[n]$ is a vertex of $\Gamma$
(equivalently, $I_\Gamma$ has no generators of degree 1).
If there are generators of degree 1 then the proof given below can
be modified by letting the index $R_1$ play the role of the index $1$. 
As $I_\Gamma$ is an almost USLI, and so
$\langle m_1, \ldots, m_l\rangle$ is a USLI,
this leaves two possible cases:

{\bf Case 1:}
{\em $m_1, \ldots, m_l$ are divisible by $x_1$, but
$m_{l+1}$ is not divisible by $x_1$.}
Since $I_\Gamma$ is squarefree strongly stable, it follows that 
$m_{l+1}=\prod_{j=2}^{d+1}x_j$. In this case each set $F_i=[n]-\{1, i\}$,
$i=2, \ldots, d+1$, is a facet of $\Gamma$. 
(Indeed the product $\prod\{x_j : j\in F_i\}$
 is not divisible by $m_{l+1}$,
and it is also not divisible by $x_1$, and hence
by $m_1, \ldots, m_{l}$, implying that $F_i$ is a face. To show that $F_i$
is a maximal face observe that 
 $F_i\cup \{i\}$ 
contains the support of $m_{l+1}$, and hence is not a face,
but then shiftedness of $\Gamma$ implies that
neither is $F_i\cup\{1\}$.)
Since there also should be a facet containing $1$, we conclude
that $\max(\Gamma)\geq d+1>\deg(m_{l+1})$, 
completing the proof of this case.

{\bf Case 2:} 
{\em All minimal generators of $I$ are divisible by $x_1$.}
In this case consider an almost USLI
$I_\Gamma':=\langle x_1, m_1/x_1, \ldots, m_{l+1}/x_1 \rangle$.
 By induction hypothesis $\Gamma'$
%(considered as a simplicial
%complex on the vertex set $[2, n]$)
has $s>\deg(m_{l+1})-1$ facets which we denote by $F_1, \ldots, F_s$.
One  easily verifies that
$\max(\Gamma)=\left\{\{1\}\cup F_1, \ldots, \{1\}\cup F_s, [2,n]\right\},
$
and so $|\max(\Gamma)|=s+1>\deg(m_{l+1})$.
\endproof

We close this section with an algebraic lemma that relates regularity of
$\Gin_{\lex}(I_\Gamma)$ to the number of facets of $\Gamma$ (for an arbitrary
complex $\Gamma$).

\begin{lemma}  \label{Pardue_lemma}
For a (finite) simplicial complex $\Gamma$, 
$\reg(\Gin_{\lex}(I_{\Gamma}))\geq |\max(\Gamma)|$.
\end{lemma}
\proof
This fact is a corollary of \cite[Lemma 23]{Pardue}
applied to squarefree (and hence radical)  ideal $I_\Gamma\in S_{[n]}$. 
For  $\sigma\subseteq[n]$, we denote by $P_\sigma$ the (prime)
ideal in $S_{[n]}$ generated by $\{x_j : j\notin\sigma\}$. It is well known
that $I_\Gamma$ has the following prime decomposition:
$
I_\Gamma=\cap_{\sigma\in\max(\Gamma)} P_\sigma.
$
 Thus the variety of $I_\Gamma$, $\mathcal{V}(I_\Gamma)$,
is the union (over $\sigma\in\max(\Gamma)$)
of the irreducible subvarieties $\mathcal{V}(P_\sigma)$.
Each such subvariety is a 
linear subspace of $\field^n$ of codimension $n-|\sigma|$.
\cite[Lemma 23]{Pardue} then implies that the monomial $m:=\prod x_i^{r_i}$,
where $r_i=|\{\sigma\in\max(\Gamma): |\sigma|=n-i\}|$,
is a minimal generator of $\Gin_{\lex}(I_\Gamma)$.
Hence $\reg(\Gin_{\lex}(I_\Gamma))\geq \deg(m)=|\max(\Gamma)|$.
\endproof

\section{Lex shifting,  $B$-numbers and the limit complex}         
                     \label{infinite_section}
In this section after defining the notion of
lexicographic shifting and the notion of $B$-numbers 
(a certain analog of the Hilbert function) we prove Theorem~\ref{Thm2}. 
We remark that extending the notion of algebraic shifting to an arbitrary term 
order
$\succ$ is not entirely automatic since the $\Phi$-image of the set of
minimal generators of $\Gin_{\succ}(I_\Gamma)\subset S_{[n]}$, 
$G(\Gin_{\succ}(I_\Gamma))$, may not be a subset of $S_{[n]}$.
 This however can be easily corrected if one considers the system of rings 
$S_{[n]}$, $n\in\N$,
endowed with natural embeddings $S_{[n]}\subseteq S_{[m]}$ 
for $m\geq n$, and makes 
all the computations in the direct limit ring 
$S=\lim_{n\rightarrow\infty}S_{[n]}=\field[x_i : i\in\N]$. This is the approach
we adopt here.
 We work with the class of monomial ideals $I\subset S$
finitely generated in each degree.
Throughout this section we use the graded 
lexicographic term order on $S$.

\begin{definition}   \label{gin_def}
Let $I$ be a monomial ideal of $S$ that is finitely generated
in each degree.
Define 
$$\Gin_{\lex}(I):=\lim_{n\rightarrow\infty}\,
\left(\Gin_{\lex}(I\cap S_{[n]})\right)S,
$$ 
where we consider $I\cap S_{[n]}$ as an ideal of $S_{[n]}$.
\end{definition}
Since the $d$-th component of $\Gin_{\lex}(I\cap S_{[n]})$ depends only on the
$d$-th component of $I\cap S_{[n]}$,
 or equivalently on the minimal generators of 
$I\cap S_{[n]}$ of degree $\leq d$,
Lemma \ref{cone}(4) implies that 
$\Gin_{\lex}(I)$ is well-defined and that for every
$d$ there is $n(d)$ such that
$(\Gin_{\lex}I)_d=((\Gin_{\lex}(I\cap S_{[n]}))S)_d$ 
for all $n\geq n(d)$.
Thus $\Gin_{\lex}(I)$ is a monomial ideal finitely generated
in each degree. (It is finitely generated if $I$ is.)
  Moreover, it follows from Lemma \ref{cone}(1) that
$\Gin_{\lex}(I)$ is a strongly stable ideal.

Recall that the squarefree operation $\Phi$ 
takes monomials of $S$ to squarefree
monomials of $S$. 
If $I\subset S$ is a monomial ideal finitely generated in each degree,
we define $\Phi(I):=\langle \Phi(m) : m\in G(I)\rangle$, where $G(I)$ is
the set of minimal generators of $I$.
\begin{definition}
Let $I$ be a homogeneous ideal of $S$ that is finitely generated
in each degree. The {\em lexicographic shifting} 
of $I$ is the squarefree strongly stable ideal 
$\Delta_{\lex}(I)=\Phi(\Gin_{\lex}(I))$. 
The {\em $i$-th lexicographic shifting} of $I$
is the ideal $\Delta_{\lex}^i(I)$, where $\Delta_{\lex}^i$ stands 
for $i$ successive applications of $\Delta_{\lex}$. 
We also define the {\em limit ideal}
$\overline{\Delta}(I):=\lim_{k\rightarrow\infty}\Delta_{\lex}^k(I)$.
\end{definition}

The rest of the section is devoted to the proof of Theorem \ref{Thm2}.
First however we digress and review
several facts on algebraic Betti numbers (defined by Eq.~(\ref{alg_betti})).

\begin{lemma}      \label{betti-prop}
Let $I$ and $J$ be monomial ideals of $S_{[n]}$.
\begin{enumerate}
\item If $I_j=J_j$ for all $0\leq j\leq j_0$, then 
$\beta_{i,j}(I)=\beta_{i,j}(J)$ for all $i$ and all $j\leq j_0$. 
\item The Betti numbers of $I\subset S_{[n]}$ coincide with those
of $IS_{[n+1]}\subset S_{[n+1]}$, that is,
 $\beta_{i,j}(I)=\beta_{i,j}(IS_{[n+1]})$ for all $i, j$.
\end{enumerate}
\end{lemma}
\proof
Part (1) follows from the standard facts that
$$\beta_{i,j}(I)=
\dim_\field \Tor_i^{S_{[n]}}(\field, I)_{j}=
\dim_\field \Tor_i^{S_{[n]}}(I, \field)_{j},$$ 
where we identify $\field$ with the $S_{[n]}$-module 
$S_{[n]}/\langle x_1, \ldots, x_n\rangle$. 
For part (2) note that if $\F$ is
the free minimal 
resolution of $I$ over $S_{[n]}$, then  $\F\otimes_{S_{[n]}} S_{[n+1]}$ 
is the free minimal resolution of $IS_{[n+1]}$ over $S_{[n+1]}$,
yielding the lemma.
% (exactness is clear and minimality
%follows from the fact that
%the maps are given by the same matrices, and hence )
\endproof

The above properties allow to extend the definition
of the Betti numbers to the class of monomial ideals of $S$
that are finitely generated in each degree.

\begin{definition}   \label{betti_def}
Let $I\subset S$ be a monomial ideal finitely generated in each
degree.
 Define 
$$\beta_{i,j}(I):=
\lim_{n\rightarrow\infty}\beta_{i,j}(I\cap S_{[n]}) \quad \mbox{for all } 
i, j\geq0,
$$
where we consider $I\cap S_{[n]}$ as an ideal of $S_{[n]}$.
\end{definition}
We remark that since $I$ is finitely generated in each degree,
for a fixed $j_0$
there exists $n_0$ such that $(I\cap S_{[n+1]})_j=((I\cap S_{[n]})S_{[n+1]})_j$
for all $0\leq j \leq j_0$ and $n\geq n_0$.
  Hence it follows from Lemma \ref{betti-prop}
that (for a fixed $i$)
the sequence $\{\beta_{i,j_0}(I\cap S_{[n]})\}_{n\in\N}$
is a constant for indices starting with $n_0$, and thus
 $\beta_{i,j_0}(I)$ is well-defined.

The Betti numbers of strongly stable ideals (of $S_{[n]}$) were computed by
Eliahou and Kervaire \cite{ElKer}, and the analog of this
formula for squarefree strongly stable ideals (of $S_{[n]}$) was established
by Aramova, Herzog, and Hibi \cite{AHHlex}. Definition \ref{betti_def}
allows to state these results as follows. (For a monomial $u$
define $m(u):=\max\{i : x_i|u\}$.)
\begin{lemma}   \label{EK}
Let $I\subset S$ be a monomial ideal finitely generated in each degree, 
let $G(I)$ denote its set of minimal generators, and let 
$G(I)_j=\{u\in G(I): \deg u=j\}$. 
\begin{enumerate}
\item If $I$ is strongly stable, then
$\beta_{i, i+j}(I)=\sum_{u\in G(I)_j} {m(u)-1 \choose i}$;
\item If $I$ is squarefree strongly stable, then
$\beta_{i, i+j}(I)=\sum_{u\in G(I)_j} {m(u)-j \choose i}$.
In particular, if $I=L(k_\bullet)$ is a USLI, then
$\beta_{i, i+j}(I)=\sum_{l=1}^{k_j}{k_1+\ldots+k_{j-1}+l-1 \choose i}$.
\end{enumerate}
\end{lemma}

Using  the notion of the Betti numbers, one can define 
a certain analog of the Hilbert function ---
 the $B$-numbers --- of a monomial ideal $I$ of $S$ that is finitely generated 
in each degree.

\begin{definition}   \label{B-definition}
Let $I\subset S$ (or $I\subset S_{[n]}$)
be a monomial ideal finitely generated in each degree, and 
let $\beta_{i,j}(I)$ be its graded Betti numbers. Define
$$
B_j(I):=\sum_{i=0}^j (-1)^i\beta_{i,j}(I) \quad \mbox{ for all }
j\geq 0 \quad  (\mbox{e.g., $B_0=0$ and $B_1(I)=|G(I)_1|$}).
$$
 The sequence 
$B(I):=\{B_j(I): j\geq 1\}$ is called the {\bf $B$-sequence} of $I$. 
\end{definition}
\begin{remark}
It is well known and is easy to prove (see \cite[Section 1B.3]{Eis2})
 that for every $n\in\N$ the polynomial
 $\sum_j B_j(I\cap S_{[n]})x^j$
equals $(1-x)^n \text{Hilb}(I\cap S_n, x)$, where 
$\text{Hilb}(I\cap S_n, x)$ is the Hilbert series of $I\cap S_{[n]}$. 
In particular, if $\Gamma$ is a $(d-1)$-dimensional
simplicial complex on $[n]$ and 
$I_{\Gamma}\subset S_{[n]}$
is its Stanley-Reisner ideal then 
\[
\frac{1-\sum_j B_j(I_\Gamma)x^j}{(1-x)^n} = \text{Hilb}(S_{[n]}/I_\Gamma, x)
=\sum_{i=0}^{d} \frac{f_{i-1}(\Gamma)x^i}{(1-x)^i}=
\frac{\sum_{i=0}^d h_i(\Gamma)x^i}{(1-x)^{d}}, 
\]
where $\{h_i(\Gamma)\}_{i=0}^d$ is the $h$-vector of $\Gamma$ \cite{St}.
(Recall that
$h_j=\sum_{i=0}^j (-1)^{j-i}{d-i \choose j-i}f_{i-1}$ for
$0\leq j \leq d$. In particular, $h_1=f_0-d$.)
Thus 
$\sum_j B_j(I_\Gamma)x^j=1-(1-x)^{h_1}\sum_i h_ix^i$
(if one assumes that $\{i\}\in\Gamma$ for every $i\in[n]$), and so
the $h$-vector of $\Gamma$ 
defines the $B$-sequence of $I_\Gamma$.
%$$\sum_j B_j(I_\Gamma S)x^j=1-(1-x)^{n-d}\sum_i h_i(\Gamma)x^i.$$
\end{remark}

The following lemma provides the
analog of the ``$f(\Gamma)=f(\Delta_{\revlex}(\Gamma))$-property".
\begin{lemma} \label{cone2}
If $I\subset S$  is a monomial ideal that is
 finitely generated in each degree, then 
the ideals $I$ and $\Delta_{\lex}(I)$ have the same $B$-sequence.
In particular, if $I$ is finitely generated, then for a sufficiently large $n$,
the ideals $I\cap S_{[n]}$ and $\Delta_{\lex}(I)\cap S_{[n]}$ 
have the same Hilbert function
(in $S_{[n]}$).
\end{lemma}
\proof
Since for every $n\in\N$ the ideals
$I\cap S_{[n]}$ and $\Gin_{\lex}(I\cap S_{[n]})$ have the same Hilbert function
(in $S_{[n]}$) (see Lemma \ref{cone}), and since 
$B_i(I)=\lim_{n\rightarrow\infty} B_i(I\cap S_{[n]})$, 
the above remark implies that $B(I)=B(\Gin_{\lex}(I))$.  Finally,
 since $\Gin_{\lex}(I)$
is a strongly stable ideal (Lemma \ref{cone}),  we infer (by comparing
the two formulas of Lemma \ref{EK})
that 
$\beta_{i,j}(\Gin_{\lex}(I))=\beta_{i,j}(\Phi\Gin_{\lex}(I))=\beta_{i,j}
(\Delta_{\lex}(I))$ for all $i, j$,
and so $B(\Gin_{\lex}(I))=B(\Delta_{\lex}(I))$.
The result follows. \endproof

Now we are ready to verify the first part of Theorem \ref{Thm2}. In fact
we prove the following slightly more general result.

\begin{theorem} \label{main}
Let $I$ be a squarefree strongly stable ideal of $S$ finitely 
generated in each degree. 
Then $\Delta_{\lex}(I)>_{\lex} I$ unless $I$ is a USLI in which case
$\Delta_{\lex}(I)= I$. 
Moreover if $I$ is finitely
generated and is not a USLI, then 
all ideals in the sequence
$\{\Delta^i_{\lex}(I)\}_{i\geq 0}$ are distinct. 
\end{theorem}

\proof There are several possible cases.

{\bf Case 1:} $I=L(k_\bullet)$ is a USLI. To prove that 
$\Delta_{\lex}(I)=I$,
it suffices to show that for every $d\geq 1$, 
$\Delta_{\lex}(L(k^{(d)})= L(k^{(d)})$, where 
$k^{(d)}:=\{k_1, \ldots, k_d, 0,0,\ldots\}$ is the sequence $k_\bullet$
truncated at $k_d$.
But this is immediate from Lemmas \ref{comb_USLI}(2) and \ref{cone2}. 
Indeed, for $n=n(d)$ sufficiently large
the simplicial complexes on the vertex set $[n]$ whose Stanley-Reisner
ideals are given by $\Delta_{\lex}(L(k^{(d)})\cap S_{[n]}$ and 
$L(k^{(d)})\cap S_{[n]}$, respectively,
are shifted and have the same $f$-numbers. 
Since the second complex is a USLI complex,
it follows that those complexes, and hence their ideals, coincide. 
 
{\bf Case 2:} $I=\langle m_1, \ldots, m_l, m_{l+1} \rangle$
 is an almost USLI.
Let $n$ be the largest index of a variable appearing in $\prod_{i=1}^{l+1}m_i$,
and let $\Gamma$ be a simplicial complex on $[n]$ with 
$I_\Gamma = I\cap S_{[n]}$.
Then
$$\reg(\Delta_{\lex}(I))=\reg(\Gin_{\lex}(I_\Gamma)) 
\stackrel{\mbox{\tiny {Lemma \ref{Pardue_lemma}}}}{\geq} |\max(\Gamma)|
\stackrel{\mbox{\tiny {Lemma \ref{main_lemma}}}}{>}\reg(I_\Gamma)=\reg(I),
$$
yielding that $\Delta_{\lex}(I)\neq I$ in this case.
Moreover, since by Eq.~(\ref{P5}),  
$\Phi(\Gin_{\revlex}(I_\Gamma))=I_\Gamma$ and since $\Phi$ is a 
lex-order preserving map,
we infer from Lemma \ref{gin>gin} that 
$\Phi(\Gin_{\lex}(I_\Gamma))\geq_{\lex} \Phi(\Gin_{\revlex}(I_\Gamma))
=I_\Gamma$,
and hence that $\Delta_{\lex}(I)>_{\lex} I$.

{\bf Case 3:} I is squarefree strongly stable, but is not a USLI. 
In this case we sort $G(I)=\{m_1, \ldots, m_l, m_{l+1}, \ldots\}$ by graded lex-order
and assume that $m_{l+1}$ is the first non-USLI generator of $I$.
Let  
$I_1=\langle m_1, \ldots, m_l \rangle$  and let 
$I_2=\langle m_1, \ldots, m_{l+1} \rangle$. 
Then $I_1$ is a USLI, $I_2$ is an almost USLI, and $I_1\subset I_2\subseteq I$. 
Hence by the previous two cases 
$I_1=\Delta_{\lex}(I_1)\subset\Delta_{\lex}(I_2)$ 
and 
$\Delta_{\lex}(I_2)>_{\lex} I_2$, and so 
 there exists a monomial $m$, $m_l>_{\lex} m>_{\lex} m_{l+1}$, such that
$m \in G(\Delta_{\lex}(I_2)) \subseteq G(\Delta_{\lex}(I))$. 
Thus
$\Delta_{\lex}(I)>_{\lex} I$. 

Finally to show that for a finitely generated ideal $I$,
all ideals in the sequence $\{\Delta^i_{\lex}(I)\}_{i\geq 0}$ are distinct, 
it suffices to check that none of those ideals is a USLI. 
This is an immediate corollary of Lemmas \ref{comb_USLI}(2) and  \ref{cone2}. \endproof

Our next goal is to prove the second part of Theorem \ref{Thm2}. 
To do that we  fix a sequence of  integers
$B=\{B_j : j\geq 1\}$ and study the class $\M(B)$ of all monomial ideals $I\subset S$ 
that are finitely generated in each degree  and satisfy $B(I)=B$.

\begin{lemma}
There is at most one USLI in the class $\M(B)$.
\end{lemma}
\proof 
Recall that  a USLI $L=L(k_\bullet)$ 
is uniquely defined by its $k$-sequence
$k_\bullet=\{k_i : i\geq 1\}$, where $k_i=\beta_{0,i}(L)=|G(L)_i|$. 
Recall also that  
$B(L)$ is a function of 
$k_\bullet$ (see Lemma \ref{EK}(2)), and so to complete the proof it suffices
to show that this function is
one-to-one, or more precisely that
 $k_j$ is determined by
$k_1, \ldots, k_{j-1}, B_j$ (for every $j\geq 1$). 
And indeed, 
%if $j=1$, then $B_1=\beta_{0,1}(J)-\beta_{1,1}(J)=\beta_{0,1}(J)=k_1$,
%and if $j>1$ 
%we have
\begin{eqnarray*}
k_j&=&\beta_{0,j}(L)=B_j-\sum_{i=1}^j (-1)^i\beta_{i,j}(L)  
   \quad (\mbox{by definition of } B_j)\\
&=& B_j-\sum_{i=1}^j(-1)^i \sum_{l=1}^{k_{j-i}}{k_1+\ldots+k_{j-i-1}+l-1 \choose i}
\quad (\mbox{by Lemma } \ref{EK}(2)).  
%&=& B_j+\sum_{i=1}^j 
%  (-1)^i\sum_{l=1}^{k_{j-i}}{k_1+\ldots+k_{j-i-1}+l-1 \choose i}
%          \quad (\mbox{by Def }\ref{USLI}).
\end{eqnarray*}
\endproof

Now we are ready to prove (the slightly more general 
version of) the second part of Theorem \ref{Thm2}.
\begin{theorem}
For every ideal $I\in\M(B)$, the limit ideal $\overline{\Delta}_{\lex}(I)$ 
is well defined and is the unique USLI of $\M(B)$.
\end{theorem}
\proof
Fix $I\in \M(B)$. To show that $\overline{\Delta}_{\lex}(I)$ 
is well defined, it suffices to check that for every $d\geq 0$,
there exists $s=s(d)$ such that 
\begin{equation}   \label{stab}
G(\Delta^{s}_{\lex}(I))_{\leq d}=G(\Delta^{s+1}_{\lex}(I))_{\leq d} 
\end{equation}
(where $G(J)_{\leq d}:=\cup_{j\leq d} G(J)_j$),
and hence that all ideals $\Delta^{i}_{\lex}(I)$, $i\geq s$,
have the same $d$-th homogeneous component.
We verify this fact by showing that the collection of all 
possible sets of minimal generators 
\begin{equation}  \label{finite}
       \mathcal{G}_{\leq d}:=\{ G(J)_{\leq d} : 
       J\in\M(B), J \mbox{ is squarefree strongly stable}\} 
 \quad \mbox{is finite}.
\end{equation}
(This yields (\ref{stab}), since all ideals $\Delta^{i}_{\lex}(I)$, $i\geq 1$,
are squarefree strongly stable, and since 
$\Delta^{i}_{\lex}(I)\leq_{\lex} \Delta^{i+1}_{\lex}(I)$ 
by Theorem \ref{main}.)
Eq.~(\ref{finite}) can be easily proved by induction.
 It clearly holds for 
$d=0$. Now if $J\in\M(B)$ is squarefree strongly stable, then
by Lemma \ref{EK}(2) and Definition \ref{B-definition}, 
$$ |G(J)_d|=\beta_{0,d}(J)=
B_d-\sum_{i=1}^{d}(-1)^i\sum_{u\in G(J)_{d-i}}{m(u)-(d-i) \choose i},
$$
so assuming that the collection 
$\mathcal{G}_{\leq d-1}$ is finite, or equivalently that the set of integers
$\{m(u): u\in G(J)_{\leq d-1}\in\mathcal{G}_{\leq d-1}\}$  is bounded 
(say by $n(d)$), 
we obtain that 
there exists a constant $g(d)$ such that
$|G(J)_d|\leq g(d)$ for all squarefree strongly stable ideals $J\in\M(B)$.
But then the squarefree strongly stable property implies that
$m(u)< n(d)+g(d)+d$ for every $u\in G(J)_{\leq d}\in \mathcal{G}_{\leq d}$,
and (\ref{finite}) follows.

The second part of the statement is now immediate:
indeed if $G(\Delta^s(I))_{\leq d} = G(\Delta^{s+1}(I))_{\leq d}$,
 then by Theorem \ref{main}, 
$G(\Delta^s(I))_{\leq d}= G(\overline{\Delta}(I))_{\leq d}$ 
is the set of minimal generators of a USLI.
\endproof

\section{Remarks on other term orders}
We close the paper by discussing several results and conjectures
related to algebraic shifting with respect to arbitrary term orders. 
 To this end, we say that an order $\succ$ 
on monomials of $S$ is a {\em term order} if 
$x_i\succ x_{i+1}$ for  $i\geq 1$, 
$m\succ m'$ as long as $\deg(m)<\deg(m')$,
and  the restriction of $\succ$
 to $S_{[n]}$
is a term order on $S_{[n]}$ for all $n\geq 1$. In addition,
we restrict our discussion
only to those term orders on $S$ that are compatible with the squarefree
operation $\Phi$, that is, $\Phi(m)\succ\Phi(m')$ if $m\succ m'$.

Similarly to Definition \ref{gin_def}, for a term order $\succ$ on $S$ and
a homogeneous ideal $I\subset S$ that is finitely generated in each degree,
we define $\Delta_\succ(I):=\Phi(\Gin_\succ(I))$. Thus $\Delta_\succ(I)$
is a squarefree strongly stable ideal that has the same $B$-sequence as $I$.
(Indeed, the proof of Lemma \ref{cone2} carries over to this more 
general case.)

We say that a squarefree monomial ideal $I\subset S$ is a {\em US$\succ$I}
if for every monomial $m\in I$ and every squarefree monomial $m'$
such that
$\deg(m)=\deg(m')$ and  $m'\succ m$, $m'$ is an element of $I$ as well. 
Being US$\succ$I implies being squarefree strongly stable.

In view of Theorems  \ref{Thm2} and \ref{AHH}
it is natural to ask the following:
\begin{enumerate}
\item Does $\Delta_\succ(I)=I$ hold for every US$\succ$I I?
\item Is there  a term order $\succ$ other than the lexicographic order
for which the equality $\Delta_\succ(I)=I$ implies that $I$ is a 
US$\succ$I?
\item Is there  a term order $\succ$ other than the 
reverse lexicographic order such that the equation $\Delta_\succ(I)=I$ 
holds for all squarefree strongly stable ideals $I$?
\end{enumerate}

The next proposition answers the first question in the affirmative.
\begin{proposition}
If $I$ is a  US$\succ$I, then $\Delta_\succ(I)=I$
for every term order on $S$ that is compatible with $\Phi$.
\end{proposition}
\proof Exactly as in the proof of Theorem \ref{main} (see the
last three lines of Case 2), 
one can show that $\Delta_\succ(I)\succeq I$. Hence either
$\Delta_\succ(I)= I$, in which case we are done, or the 
$\succ$-largest monomial, $m$, in the symmetric difference of 
$G(\Delta_\succ(I))$ and $G(I)$ is an element of $G(\Delta_\succ(I))$.
Since $I$ is a US$\succ$I, we obtain in the latter case that 
$G(\Delta_\succ(I))_i=G(I)_i$ for all $i<\deg(m)$ and
$$
G(I)_{i_0}=\{m'\in G(\Delta_\succ(I))_{i_0} : m'\succ m\}
%\subseteq G(\Delta_\succ(I))_{i_0}-\{m\} 
\quad \mbox{ for } i_0=\deg(m),
$$
that is, $G(I)_{i_0}$ is a strict subset of $ G(\Delta_\succ(I))_{i_0}$.
This is however impossible, since it contradicts the fact that
the ideals $I$ and $\Delta_\succ(I)$ have the same $B$-sequence.
\endproof

The answer to the second question is negative as follows from 
the following result.

\begin{proposition}
If $I$ is a USLI, then $\Delta_\succ(I)=I$ for all term orders $\succ$.
\end{proposition}
We omit the proof as it is completely analogous to that of 
Theorem \ref{main}, Case 1.

While we do not know the answer to the third question, we believe that it
is negative. In fact it is tempting to conjecture that the following holds.
Let $\succ$ be a term order on $S$ other than the (graded) 
reverse lexicographic order, and let $k\geq 2$ be the smallest degree
on which $\succ$ and revlex disagree. Write $m_i$ to denote the $i$th 
squarefree monomial of $S$ of degree $k$ with respect to the revlex order.
(It is a fundamental property of the revlex order that every squarefree 
monomial of $S$ of degree $k$ is of the form $m_i$ for some finite $i$.)

\begin{conjecture}
Let $i_0\geq 1$ be the smallest index for which 
$I_{i_0}:=\langle m_1, \cdots, m_{i_0}\rangle$ is not a US$\succ$I. 
Then $\Delta_{\succ}(I_{i_0})\neq I_{i_0}$.
\end{conjecture}

\section*{Acknowledgments} We are grateful to Aldo Conca for
helpful discussions and to the anonymous referees for insightful comments.

\end{document}